\documentclass[12pt]{amsart}

\begin{document}
\setcounter{page}{1}

\title{Riemann's hypothesis via Robin's Theorem}
\author[D. Martila]{Dmitri Martila}
\author[S. Groote]{Stefan Groote}
\address{Institute of Physics, University of Tartu\\
W.~Ostwaldi~1, EE-50411 Tartu, Estonia}
\thanks{eestidima@gmail.com, Independent Researcher\\
J. V. Jannseni 6--7, P\"arnu 80032, Estonia}

\begin{abstract}
As essential condition for the validy of Robin's Theorem as a precondition
for the proof of the Riemann hypothesis, we show that the
minimum of the function $F={\rm e}^{\gamma}\,\ln(\ln\,n)-\sigma(n)/n$ is
found to be positive. Therefore, $F>0$ holds for any $n>5040$.

MSC Class: 11M26, 11M06.
\end{abstract}

\maketitle

Robin's Theorem~\cite{Robin} states that if
\begin{equation}\label{F}
F={\rm e}^{\gamma}\,\ln(\ln\,n)-\frac{\sigma(n)}{n}>0
\end{equation}
for $n> 5040$, where $\gamma\approx 0.577$ is the Euler--Mascheroni constant
and $\sigma(n)$ is the sum-of-divisors function, the Riemann hypothesis is true.

In the following we use methods of functional analysis to show that the minimum
of $F$ has to be positive. According to the fundamental theorem of arithmetic,
one has 
\begin{equation}\label{n}
n=\prod_{i=1}^{\kappa} p_i^{\alpha_i}\,. 
\end{equation}
Hence, $F$ is a unique function of the powers $\alpha_g$ via
$n=n_0\,p_g^{\alpha_g}$, where $p_g$ is not a divisor of $n_0$.
Methodologically, derivatives with respect to $\alpha_g$ can be taken by
calculating finite differences~\cite{Fin}. Therefore,
\begin{equation}
\frac{\Delta p_g^{\alpha_g}}{\Delta\alpha_g}=
p_g^{\alpha_g}-p_g^{\alpha_g-1}=p_g^{\alpha_g}\,(1-1/p_g)
\end{equation}
and, accordingly,
\begin{equation}
S(n)=\frac{\Delta}{\Delta\alpha_g}\,\ln(\ln\,n)=
\frac{1}{n\,\ln\,n}\,\frac{\Delta n}{\Delta\alpha_g}=\frac{1-1/p_g}{\ln\,n}\,. 
\end{equation}
On the other hand,
\begin{equation}
\frac{\sigma(n)}n=A\,u(\alpha_g)\,,\qquad
u(\alpha_g)=1+\frac{1}{p_g}+\frac{1}{p_g^2}+...\frac{1}{p_g^{\alpha_g}}\,. 
\end{equation}
$p_g$ is not a divisor of $A$, as it is made obvious by a simple example in the
Appendix. Again, the finite difference results in
\begin{equation}
\frac\Delta{\Delta\alpha_g}\left(\frac{\sigma(n)}n\right)
  =A\,(u(\alpha_g)-u(\alpha_g-1))=\frac{A}{p_g^{\alpha_g}}\,.
\end{equation}
By analysis, the minimum of $F$ can be found via 
\begin{equation}
0=\frac{\Delta F}{\Delta\alpha_g}
  ={\rm e}^{\gamma}\,\frac{1-1/p_g}{\ln\,n}-\frac{A}{p_g^{\alpha_g}}
  ={\rm e}^{\gamma}\,\frac{1-1/p_g}{\ln\,n}
  -\frac{\sigma(n)}{n\,u(\alpha_g)\,p_g^{\alpha_g}}=0\,.
\end{equation}
This means
\begin{equation}
{\rm e}^{\gamma}\,(1-1/p_g)\,(1+p_g+p_g^2+...+p_g^{\alpha_g})
  =\sigma(n)\,\frac{\ln\,n}{n}\,,
\label{FF34}
\end{equation}
and resummed
\begin{equation}
{\rm e}^{\gamma}\,(1-1/p_g)\,\frac{p_g^{\alpha_g+1}-1}{p_g-1}
  ={\rm e}^{\gamma}\left(p_g^{\alpha_g}-\frac{1}{p_g}\right)
  =\sigma(n)\,\frac{\ln\,n}{n}
\end{equation}
or
\begin{equation}
p_g^{\alpha_g}
  ={\rm e}^{-\gamma}\,\sigma(n)\,\frac{\ln\,n}{n}+\frac{1}{p_g}\,.
\label{gty}
\end{equation}
For a given $n$, this formula provides values for $p_g$ and $\alpha_g$. The
higher $p_g$ is chosen, the lower $\alpha_g$ is found. The highest prime
$g=\kappa$ has $\alpha_\kappa=1$. From Eq.~(\ref{FF34}), one has
\begin{equation}
\frac{\sigma(n)}n={\rm e}^{\gamma}\,
  \frac{p_{\kappa}^2-1}{p_{\kappa}\,\ln\,n}\,,
\end{equation}
and inserting to Eq.~(\ref{F}) one obtains
 \begin{equation}
F\,p_{\kappa}\,\ln\,n={\rm e}^{\gamma}\,p_{\kappa}\,\ln\,n\,\ln(\ln\,n)
  -{\rm e}^{\gamma}\,(p_{\kappa}^2-1)\,,
\label{hjk}
\end{equation}
where $n$ is a potential counter-example. As the powers in Eq.~(\ref{n}) up to
$i=\kappa$ are at least $1$, $\alpha_i\ge 1$, Eq.~(\ref{gty}) gives
\begin{equation}
n\ge\prod_{i=1}^{\kappa}p_i=\exp(\theta(p_{\kappa}))\,,
\end{equation}
where $\theta$ is the first Chebyshev function with
$\lim_{x\to\infty}(\theta(x)/x)=1$~\cite{Tom}. Therefore,
$n>\exp(p_{\kappa}/2)$ inequality, and the right hand of the Eq.~(\ref{hjk})
is positive. This means, that $F$ cannot have negative values.

\begin{appendix}
\section*{Appendix: on the sum-of-divisors function}
\setcounter{equation}{0}\def\theequation{A\arabic{equation}}
As an example, we consider $n=28=2\cdot 14=2\cdot 2\cdot 7=2^2\cdot 7$.
Therefore, $\sigma(28)=(1+7)(1+2+2^2)=1+2+4+7+14+28$ and
\begin{equation}
\frac{\sigma(28)}{28}=\frac{1+7}7\cdot\frac{1+2+2^2}{2^2}\,.
\end{equation}
Selecting $p_g=2$, one has $A=(1+7)/7$ and $u_g=(1+2+2^2)/2^2$.
\end{appendix}

\end{document}